\documentclass[12pt]{amsart}
\usepackage{amsmath,amssymb, euscript, graphicx}

\theoremstyle{definition}

\begin{document}

 \address{Azer Akhmedov, Department of Mathematics,
North Dakota State University,
Fargo, ND, 58108, USA}
\email{azer.akhmedov@ndsu.edu}

  \begin{center} {\Large \bf Groups not acting on compact metric spaces by homeomorphisms} \end{center}
  
  \medskip
  
  \begin{center} Azer Akhmedov \end{center}
  
  \medskip
  
  {\bf Abstract:} {\Small We show that the direct sum of uncountably many non-Abelian groups does not embed into the group of homeomorphisms of a compact metric space.} 
   
\vspace{1cm}

  The main result of this note does not seem to exist in the literature. Among other things, it provides very simple examples of left-orderable groups of continuum cardinality which do not embed in $\mathrm{Homeo}_{+}(\mathbb{R})$ \footnote{It is well known that a countable group is left-orderable iff it embeds into $\mathrm{Homeo}_{+}(\mathbb{R})$. Thus it becomes an interesting question whether or not there exists a left-orderable group of continuum cardinality which does not embed in $\mathrm{Homeo}_{+}(\mathbb{R})$. This question has been addressed in a very recent paper [M] where two different examples have been constructed.}

 \medskip
 
  {\bf Theorem 1.} Let $X$ be a compact metric space, $I$ be an uncountable set, and $G_{\alpha }$ be a non-Abelian group, for all $\alpha \in I$. Then the direct sum $\displaystyle \mathop{\oplus }_{\alpha \in I}G_{\alpha }$ does not embed into the group of homeomorphisms of $X$.
  
\medskip

  {\bf Proof.} We will assume that the direct sum $\displaystyle \mathop{\oplus }_{\alpha \in I}G_{\alpha }$ embeds in $\mathrm{Homeo}(X)$. 
  
  \medskip
  
  For all $\alpha \in I$, let $f_{\alpha }, g_{\alpha }$ be some non-commuting elements in $G_{\alpha }$. We will denote the metric on $X$ by $d(.,.)$. Notice that the group $\mathrm{Homeo}(X)$ of homeomorphisms of $X$ also becomes a metric space where for $\phi , \psi \in \mathrm{Homeo}(X)$, we define the metric by $D(\phi ,\psi ) = \displaystyle \mathop{\sup }_{x\in X}d(\phi (x), \psi (x))$.
  
  \medskip 
  
  By compactness of $X$, we have a sequence  $C_n, n\geq 1$ of finite subsets of $X$ such that 
 
  (i) for all $n\geq 1$, $C_n$ is a $\frac{1}{n}$-net, i.e. for all $x\in X$, there exists $z\in C_n$ such that $d(z,x) < \frac{1}{n}$, 
   
  (ii) $C_n\subseteq C_{n+1}$, for all $n\geq 1$, 
  
  (iii) $\overline {{\cup }_{n\geq 1}C_n} = X$.
  
  \medskip
  
  Let $C_n = \{x_1^{(n)}, \dots , x_{p_n}^{(n)}\}, n\geq 1$. For all $n\geq 1$ and $\delta > 0$, we will also write $$\mathcal{F}_{n,\delta } = \{\phi \in \mathrm{Homeo}(X) \ | \ \mathrm{diam}\phi (B_{\frac{2}{n}}(x_i^{(n)})) < \delta , \forall i\in \{1, \dots , p_n\}\}$$   
  
  \medskip
  
  Since $\mathrm{Homeo}(X)\backslash \{1\} = \displaystyle \mathop{\cup }_{i\geq 1}\{\phi \ | \ D(\phi ,1) > \frac{1}{i}\}$ there exists $c > 0$ and an uncountable set $I_1\subset I$ such that for all $\alpha \in I_1$ we have  $D([f_{\alpha },g_{\alpha }], 1) > c.$ 
  
  \medskip
  
  By compactness, all homeomorphisms of $X$ are uniformly continuous. Then there exists an uncountable set $I_2\subset I_1$ and positive integers $n, m$ such that $\max \{\frac{1}{n}, \frac{1}{m}\} < \frac{c}{10}$ and for all $\beta \in I_2$ we have $\{f_{\beta }, g_{\beta }\}\subseteq \mathcal{F}_{n , \frac{1}{m}}$.
  
  \medskip
  
  Then, since $\mathrm{Homeo}(X)$ is separable, there exists an uncountable set $I_3\subset I_2$ and $f_{\star }, g_{\star }\in \mathrm{Homeo}(X)$ such that   for all $\gamma \in I_3$ we have $\max \{D(f_{\gamma }, f_{\star }), D(g_{\gamma }, g_{\star })\} < \frac{1}{100n}$.
  
  \medskip
   
  Then for all $\gamma , \eta \in I_3$ and for all $x\in X$ we have the inequalities $$d(g_{\eta }f_{\gamma }(x), g_{\gamma }f_{\gamma }(x)) < \frac{1}{50n} \ \mathrm{and} \  d(f_{\gamma }g_{\eta }(x), f_{\gamma }g_{\gamma }(x)) < \frac{c}{5}.$$ On the other hand, for any two distinct $\eta , \gamma \in I_3$ there exists $x_0\in X$ such that $d(f_{\gamma }g_{\gamma }(x_0), g_{\gamma }f_{\gamma }(x_0)) > c$. 
  
  \medskip
  
  By triangle inequality, we obtain that $d(g_{\eta }f_{\gamma }(x_0), f_{\gamma }g_{\eta }(x_0)) > \frac{c}{2}$. Thus $[g_{\eta }, f_{\gamma }]\neq 1$. Contradiction. $\square $
  
  \medskip
  
  {\bf Remark.} To obtain interesting applications of the main theorem, $X$ can be taken to be an arbitrary compact manifold (with boundary). In the case of a closed interval $I\cong [0,1]$, we obtain a non-embeddability result into $\mathrm{Homeo}_{+}(\mathbb{R})$, i.e. the direct sum of uncountably many isomorphic copies of a non-Abelian group $G$ does not embed in $\mathrm{Homeo}_{+}(\mathbb{R})$. The result of Theorem 1 can be extended to a much larger class of topological spaces. One can also drop or replace the non-Abelian condition on the group $G$ in certain special contexts. 
  
    \vspace{1cm}
    
    {\bf References:}
    
    [M] K.Mann, Left-orderable groups that don't act on the line. {\em Mathematische Zeitschrift} {\bf vol. 280}, issue 3, (2015) 905-918.  
  
 \end{document}